\documentclass[10pt,twoside]{article}
\usepackage{graphicx}
\usepackage{amsmath}
\usepackage{Latex-document}

\title{\bf Popularizing Mathematics: \vskip -2mm
From Eight to Infinity \vskip 6mm}

\author{V. L. Hansen\vspace*{-0.5cm}\thanks{Department of Mathematics,
Technical University of Denmark, Matematiktorvet, Bygning 303,
DK-2800 Kgs. Lyngby, Denmark. E-mail: V.L.Hansen@mat.dtu.dk}}

\date{\vspace{-8mm}}

\begin{document}

\maketitle

\thispagestyle{first} \setcounter{page}{885}

\begin{abstract}

\vskip 3mm

It is rare to succeed in getting mathematics into ordinary
conversation without meeting all kinds of reservations. In order
to raise public awareness of mathematics effectively, it is
necessary to modify such attitudes. In this paper, we point to
some possible topics for general mathematical conversation.

\vskip 4.5mm


\end{abstract}

\vskip 12mm

\section{Introduction}\label{section 1}\setzero

\vskip-5mm \hspace{5mm}

``One side will make you grow taller, and the other side will make you grow shorter", says the caterpillar to
Alice in the fantasy {\em Alice's Adventures in Wonderland} as it gets down and crawls away from a marvellous
mushroom upon which it has been sitting. In an earlier episode, Alice has been reduced to a height of three inches
and she would like to regain her height. Therefore she breaks off a piece of each of the two sides of the
mushroom, incidentally something she has difficulties identifying, since the mushroom is entirely round. First
Alice takes a bit of one of the pieces and gets a shock when her chin slams down on her foot. In a hurry she eats
a bit of the other piece and shoots up to become taller than the trees. Alice now discovers that she can get
exactly the height she desires by carefully eating soon from one piece of the mushroom and soon from the other,
alternating between getting taller and shorter, and finally regaining her normal height.

Few people link this to mathematics, but it reflects a result
obtained in 1837 by the German mathematician Dirichlet, on what is
nowadays known as {\em conditionally convergent infinite series}
--- namely, that one can assign an arbitrary value to infinite sums
with alternating signs of magnitudes tending to zero, by changing
the order in which the magnitudes are added. Yes, mathematics can
indeed be fanciful, and as a matter of fact, {\em Alice's
Adventures in Wonderland} was written, under the pseudonym Lewis
Carroll, by an English mathematician at the University of Oxford
in 1865. The fantasy about Alice is not the only place in
literature where you can find mathematics. And recently,
mathematics has even entered into stage plays \cite{plays} and
movies \cite{movies}.

\section{On the special position of mathematics} \label{section 2}\setzero

\vskip-5mm \hspace{5mm}

When the opportunity arises, it can be fruitful to incorporate extracts from literature or examples from the arts
in the teaching and dissemination of mathematics. Only in this way can mathematics eventually find a place in
relaxed conversations among laymen without immediately being rejected as incomprehensible and relegated to
strictly mathematical social contexts.

Mathematics occupies a special position among the sciences and in
the educational system. This position is determined by the fact
that mathematics is an {\em a priori} science building on ideal
elements abstracted from sensory experiences, and at the same time
mathematics is intimately connected to the experimental sciences,
traditionally not least the natural sciences and the engineering
sciences. Mathematics can be decisive when formulating theories
giving insight into observed phenomena, and often forms the basis
for further conquests in these sciences because of its power for
deduction and calculation. The revolution in the natural sciences
in the 1600s and the subsequent technological conquests were to an
overwhelming degree based on mathematics. The unsurpassed strength
of mathematics in the description of phenomena from the outside
world lies in the fascinating interplay between the concrete and
the abstract.

In the teaching of mathematics, and when explaining the essence of
mathematics to the public, it is important to get the abstract
structures in mathematics linked to concrete manifestations of
mathematical relations in the outside world. Maybe the impression
can then be avoided that abstraction in mathematics is falsely
identified with pure mathematics, and concretization in
mathematics just as falsely with applied mathematics. The booklet
\cite{mireille} is a contribution in that spirit.

\section{On the element of surprise in mathematics}
\label{section 3} \setzero

\vskip-5mm \hspace{5mm}

Without a doubt, mathematics makes the longest-lasting impression when it is used to explain a counter-intuitive
phenomenon. The element of surprise in mathematics therefore deserves particular attention.

As an example, most people --- even teachers of mathematics ---
find it close to unbelievable that a rope around the Earth along
the equator has to be only $2\pi$ metres (i.e., a little more than
6 metres) longer, in order to hang in the air 1 metre over the
surface all the way around the globe.

Another easy example of a surprising fact in mathematics can be
found in connection with figures of constant width. At first one
probably thinks that this property is restricted to the circle.
But if you round off an equilateral triangle by substituting each
of its sides with the circular arc centred at the opposite corner,
one obtains a figure with constant width. This rounded triangle is
called a {\em Reuleaux triangle} after its inventor, the German
engineer Reuleaux. It has had several technological applications
and was among others exploited by the German engineer Wankel in
his construction of an internal combustion engine in 1957.
Corresponding figures with constant width can be constructed from
regular polygons with an arbitrary odd number of edges. In the
United Kingdom, the regular heptagon has been the starting point
for rounded heptagonal coins (20p and 50p).

Figures of constant width were also used by the physicist Richard
Feynmann to illustrate the dangers in adapting figures to given
measurements without knowledge of the shape of the figure.
Feynmann was a member of the commission appointed to investigate
the possible causes that the space shuttle Challenger on its tenth
mission, on 28 January 1986, exploded shortly after take off. A
problematic adaption of figures might have caused a leaky assembly
in one of the lifting rockets.

\section{Mathematics in symbols} \label{section 4} \setzero

\vskip-5mm \hspace{5mm}

Mathematics in symbols is a topic offering
good possibilities for conversations with a mathematical touch.
Through symbols, mathematics may serve as a common language by
which you can convey a message in a world with many different
ordinary languages. But only rarely do you find any mention of
mathematics in art catalogues or in other contexts where the
symbols appear. The octagon, as a symbol for eternity, is only one
of many examples of this.

\subsection*{The eternal octagon}\setzero

\vskip-5mm \hspace{5mm}

Once you notice it, you will find the octagon in very many places
- in domes, cupolas and spires, often in religious sanctuaries.
This representation is strong in Castel del Monte, built in the
Bari province in Italy in the first half of the 1200s for the holy
Roman emperor Friedrich II. The castle has the shape of an
octagon, and at each of the eight corners there is an octagonal
tower \cite{gotze}. The octagon is also found in the beautiful
Al-Aqsa Mosque in Jerusalem from around the year 700, considered
to be one of the three most important sites in Islam.

In Christian art and architecture, the octagon is a symbol for the
{\em eighth day} (in Latin, {\sl octave dies}), which gives the
day, when the risen Christ appeared to his disciples for the
second time after his resurrection on Easter Sunday. In the Jewish
counting of the week, in use at the time of Christ, Sunday is the
first day of the week, and hence the eighth day is the Sunday
after Easter Sunday. In a much favoured interpretation by the
Catholic Fathers of the Church, the Christian Sunday is both the
first day of the week and the eighth day of the week, and every
Sunday is a celebration of the resurrection of Christ, which is
combined with the hope of eternal life. Also in Islam, the number
eight is a symbol of eternal life.

Such interpretations are rooted in old traditions associated with
the number eight in oriental and antique beliefs. According to
Babylonian beliefs, the soul wanders after physical death through
the seven heavens, which corresponds to the material heavenly
bodies in the Babylonian picture of the universe, eventually to
reach the eighth and highest heaven. In Christian adaptation of
these traditions, the number eight (octave) becomes the symbol for
the eternal salvation and fuses with eternity, or infinity. This
is reflected in the mathematical symbol for infinity, a figure
eight lying down, used for the first time in 1655 by the English
mathematician John Wallis.

\subsection*{The regular polyhedra}\setzero

\vskip-5mm \hspace{5mm}

Polyhedra fascinate many people. A delightful and comprehensive
study of polyhedra has been made by Cromwell \cite{crom}.

Already Pythagoras in the 500s BC knew that there can exist only
five types of regular polyhedra, with respectively 4 (tetra), 6
(hexa), 8 (octa), 12 (dodeca) and 20 (icosa) polygonal lateral
faces. In the dialogue {\em Timaeus}, Plato associated the {\em
tetrahedron}, the {\em octahedron}, the {\em hexahedron} [cube]
and the {\em icosahedron} with the four elements in Greek
philosophy, fire, air, earth, and water, respectively, while in
the {\em dodecahedron}, he saw an image of the universe itself.
This has inspired some globe makers to represent the universe in
the shape of a dodecahedron.

It is not difficult to speculate as to whether there exist more
than the five regular polyhedra. Neither is it difficult to
convince people that no kinds of experiments are sufficient to get
a final answer to this question. It can only be settled by a
mathematical proof. A proof can be based on the theorem that the
alternating sum of the number of vertices, edges and polygonal
faces in the surface of a convex polyhedron equals $2$, stated by
Euler in 1750; see e.g. \cite{hans1}.

\subsection*{The cuboctahedron}\setzero

\vskip-5mm \hspace{5mm}

There are several other polyhedra that are ascribed a symbolic
meaning in various cultures. A single example has to suffice.

If the eight vertices are cut off a cube by planes through the
midpoints in the twelve edges in the cube one gets a polyhedron
with eight equilateral triangles and six squares as lateral faces.
This polyhedron, which is easy to construct, is known as the {\em
cuboctahedron}, since dually it can also be constructed from an
octahedron by cutting off the six vertices in the octahedron in a
similar manner. The cuboctahedron is one of the thirteen
semiregular polyhedra that have been known since Archimedes. In
Japan cuboctahedra have been widely used as decorations in
furniture and buildings. Lamps in the shape of cuboctahedra were
used in Japan already in the 1200s, and they are still used today
in certain religious ceremonies in memory of the dead
\cite{miyazaki}.

\subsection*{The yin-yang symbol}\setzero

\vskip-5mm \hspace{5mm}

Yin and yang are old principles in Chinese cosmology and
philosophy that represent the dark and the light, night and day,
female and male. Originally yin indicated a northern hill side,
where the sun does not shine, while yang is the south side of the
hill. Everything in the world is viewed as an interaction between
yin and yang, as found among others in the famous oracular book
{\em I Ching} ({\em Book of Changes}), central in the teaching of
Confucius (551--479 BC). In this ancient Chinese system of
divination, an oracle can be cast by flipping three coins and the
oracle is one of sixty-four different hexagrams each composed of
two trigrams. The three lines in a trigram are either straight
(yang) or broken (yin), thus giving eight different trigrams and
sixty-four hexagrams. In Taoism, the eight trigrams are linked to
immortality. Yin and Yang as philosophical notions date back to
the 400s BC.

The mathematical symbol for yin-yang is a circle divided into two
equal parts by a curve made up of two smaller semicircles with
their centres on a diameter of the larger circle. The mathematics
in this beautiful and very well known symbol is simple, but
nevertheless it contains some basic geometrical forms, and thereby
offers possibilities for mathematical conversations.

\subsection*{The Borromean link}\setzero

\vskip-5mm \hspace{5mm}

Braids, knots and links form a good topic
for raising public awareness of mathematics, as demonstrated in
the CD-Rom \cite{brown}.

Here we shall only mention the Borromean link. This fascinating
link consists of a system of three interlocking rings, in which no
pair of rings interlocks. In other words, the three rings in a
Borromean link completely falls apart if any one of the rings is
removed from the system.

The Borromean link is named after the Italian noble family
Borromeo, who gained a fortune by trade and banking in Milano in
the beginning of the 1400s. In the link found in the coat of arms
of the Borromeo family, the three rings apparently interlock in
the way described, but a careful study reveals that the link often
has been changed so that it does not completely fall apart if an
arbitrary one of the rings is removed from the system. The link in
the coat of arms of the Borromeo family is a symbol of
collaboration.

\section{Mathematics in art and design}
\label{section 5} \setzero

\vskip-5mm \hspace{5mm}

Artists, architects and designers give life
to abstract ideas in concrete works of art, buildings, furniture,
jewellery, tools for daily life, or, by presenting human
activities and phenomena from the real world, dynamically in films
and television. The visual expressions are realized in the form of
paintings, images, sculptures, etc. From a mathematical point of
view, all these expressions represent geometrical figures.

For a mathematician, the emphasis is on finding the abstract forms
behind the concrete figures. In contrast, the emphasis of an
artist or a designer is to realize the abstract forms in concrete
figures. At the philosophical level many interesting conversations
and discussions about mathematics and its manifestations in the
visual arts can take place following this line of thought.

More down to earth, the mathematics of perspective invites many
explanations of concrete works of arts \cite{field}, and the
beautiful patterns in Islamic art inspire discussions on geometry
and symmetry \cite{abas}. In relation to the works of the Dutch
artist Escher, it is possible to enter into conversations about
mathematics at a relatively advanced level, such as the Poincar\'e
disc model of the hyperbolic plane \cite{emmer}, \cite{escher}.
The sculptures of the Australian-British artist John Robinson ---
such as his wonderful sculpture {\em Immortality} that exhibits a
M{\"o}bius band shaped as a clover-leaf knot - offer good
possibilities for conversations on knots \cite{brown}.

In his fascinating book \cite{wilson}, Wilson tells many stories
of mathematics in connection with mathematical motifs on postage
stamps.

\section{Mathematics in nature} \label{section 6} \setzero

\vskip-5mm \hspace{5mm}

The old Greek thinkers --- in particular Plato (427--347 BC) ---
were convinced that nature follows mathematical laws. This belief
has dominated thinking about natural phenomena ever since, as
witnessed so strongly in major works by physicists, such as
Galileo Galilei in the early 1600s, Isaac Newton in his work on
gravitation in {\em Principia Mathematica} 1687, and Maxwell in
his major work on electromagnetism in 1865; cf. \cite{hans1}. In
his remarkable book {\em On Growth and Form} (1917), the zoologist
D'Arcy Wentworth Thompson concluded that wherever we cast our
glance we find beautiful geometrical forms in nature. With this
point of departure, many conversations on mathematics at different
levels are possible. I shall discuss a few phenomena that can be
described in broad terms without assuming special mathematical
knowledge, but nevertheless point to fairly advanced mathematics.

\subsection*{Spiral curves and the spider's web}\setzero

\vskip-5mm \hspace{5mm}

Two basic motions of a point in the
Euclidean plane are motion in a line, and periodic motion
(rotation) around a central point. Combining these motions, one
gets spiral motions along {\em spiral curves} around a {\em spiral
point.} Constant velocity in both the linear motion and the
rotation gives an {\em Archimedean spiral} around the spiral
point. Linear motion with exponentially growing velocity and
rotation with constant velocity gives a {\em logarithmic spiral}.

A logarithmic spiral is also known as an {\em equiangular spiral},
since it has the following characteristic property: {\em the angle
between the tangent to the spiral and the line to the spiral point
is constant.}

Approximate Archimedean and logarithmic spirals enter in the
spider's construction of its web. First, the spider constructs a
Y-shaped figure of threads fastened to fixed positions in the
surroundings, and meeting at the centre of the web. Next, it
constructs a frame around the centre of the web and then a system
of radial threads to the frame of the web, temporarily held
together by a non-sticky logarithmic spiral, which it constructs
by working its way out from the centre to the frame of the web.
Finally, the spider works its way back to the centre along a
sticky Archimedean spiral, while eating the logarithmic spiral
used during the construction.

\subsection*{Helices, twining plants and an optical illusion}\setzero

\vskip-5mm \hspace{5mm}

A space curve on a cylinder, for which all
the tangent lines intersects the generators of the cylinder in a
constant angle, is a {\em helix}. This is the curve followed by
twining plants such as bindweed (right-handed helix) and hops
(left-handed helix).

Looking up a long circular cylinder in the wall of a circular
tower enclosing a spiral staircase, all the generators of the
cylinder appear to form a system of lines in the ceiling of the
tower radiating from the central point, and the helix in the
banister of the spiral staircase appears as a logarithmic
(equiangular) spiral. A magnificent instance of this is found in
the spiral staircase in Museo do Popo Galergo, Santiago, Spain,
where you have three spiral staircases in the same tower.

\subsection*{Curvature and growth phenomena in nature}\setzero

\vskip-5mm \hspace{5mm}

Many growth phenomena in nature exhibit
curvature. As an example, we discuss the shape of the shell of the
primitive cuttlefish nautilus.

The fastest way to introduce the notion of curvature at a
particular point of a plane curve is by approximating the curve as
closely as possible in a neighbourhood of the point with a circle,
called the {\em circle of curvature}, and then measuring the
curvature as the reciprocal of the radius in the approximating
circle. If the curve is flat in a neighborhood of the point, or if
the point is an inflection point, the approximating circle
degenerates to a straight line, the {\em tangent} of the curve; in
such situations, the curvature is set to $0$. Finally, the
curvature of a curve has a sign: positive if the curve turns to
the left (anticlockwise) in a neighbourhood of the point in
question, and negative if it turns to the right (clockwise).

Now to the mathematics of a nautilus shell. We start out with an
equiangular spiral. At each point of the spiral, take the circle
of curvature and replace it with the similar circle centred at the
spiral and orthogonal to both the spiral and the plane of the
spiral. For exactly one angle in the equiangular spiral, the
resulting surface winds up into a solid with the outside of one
layer exactly fitting against the inside of the next. The surface
looks like the shell of a nautilus.

In the Spring of 2001, I discovered that this unique angle is very
close to half the golden angle (the smaller of the two angles
obtained by dividing the circumference of the unit circle in
golden ratio), so close in fact that it just allows for thickness
of the shell. More details can be found on my home page
\cite{hans5}.

\section{Geometry of space} \label{section 7} \setzero

\vskip-5mm \hspace{5mm}

Questions about geometries alternative to
Euclidean geometry easily arise in philosophical discussions on
the nature of space and in popularizations of physics. It is
difficult material to disseminate, but still, one can come far by
immediately incorporating a concrete model of a non-Euclidean
geometry into the considerations. In this respect, I find the disc
model of {\em the hyperbolic plane} suggested by Poincar\'e (1887)
particularly useful. In addition to its mathematical uses, it was
the inspiration for Escher in his four woodcuts Circle Limit
I--IV, (\cite{escher}, page 180).

In (\cite{hans4}, Chapter 4), I gave a short account of the
Poincar\'e disc model, which has successfully been introduced in
Denmark at the upper high-school level, and where you encounter
some of the surprising relations in hyperbolic geometry. In
particular, the striking difference between tilings with congruent
regular $n$-gons in the Euclidean plane, where you can tile with
such $n$-gons only for $n = 3, 4, 6$, and in the hyperbolic plane,
where you can tile with such $n$-gons for each integer $n \ge 3$.

At a more advanced level, one can introduce a hyperbolic structure
on surfaces topologically equivalent to the surface of a sphere
with $p$ handles, for each {\em genus} $p \ge 2$, by pairwise
identification of the edges in a regular hyperbolic $4p$-gon; cf.
\cite{hans2}.

The appearance of non-Euclidean geometries raised the question of
which geometry provides the best model of the physical world. This
question was illustrated in an inspired poster designed by Nadja
Kutz; cf. \cite{emsposters}.

\section{Eternity and infinity}\label{section 8} \setzero

\vskip-5mm \hspace{5mm}

Eternity and, more generally, infinity are
notions both expressing the absence of limits. Starting from the
classical paradoxes of Zeno, which were designed to show that the
sensory world is an illusion, there are good possibilities for
conversations about mathematical notions associated with infinity.
And one can get far even with fairly difficult topics such as
infinite sums; cf. \cite{hans6}.

\subsection*{The harmonic series}
\hspace*{8mm} A conversation about infinite sums (series)
inevitably gets on to the {\em harmonic} series:
$$
\sum_{n=1}^{\infty}{\frac1n} = 1 + \frac 12 + \frac 13 + \frac 14 + \dots + \frac 1n + \cdots.
$$
As you know, this series is {\em divergent}, i.e. the $N$th {\em
partial sum}
$$
S_N = 1 + \frac 12 + \frac 13 + \frac 14 + \dots + \frac 1N
$$
increases beyond any bound, for increasing $N$.

The following proof of this fact made a very deep impression on me
when I first met it. First notice, that for each number $k$, the
part of the series from $\frac 1{k+1}$ to $\frac 1{2k}$  contains
$k$ terms, each greater than or equal to $\frac 1{2k}$, so that
$$
\frac 1{k+1} + \frac 1{k+2} + \dots + \frac 1{2k} \  \ge \ \frac
1{2k} + \frac 1{2k} + \dots + \frac 1{2k} = \frac 12.
$$
With this fact at your disposal, it is not difficult to divide the
harmonic series into infinitely many parts, each greater than or
equal to $\frac 12$, proving that the partial sums in the series
grow beyond any bound when more and more terms are added.

\subsection*{The alternating harmonic series}\setzero

\vskip-5mm \hspace{5mm}

The {\em alternating} harmonic series
$$
\sum_{n=1}^{\infty}{(-1)^{n-1}\frac1n} = 1 - \frac 12 + \frac 13 -
\frac 14 + \dots + (-1)^{n-1} \frac 1n + \dots $$ offers great
surprises. As shown in 1837 by the German mathematician Dirichlet,
this series can be {\em rearranged}, that is, the order of the
terms can be changed, so that the rearrangement is a convergent
series with a sum arbitrarily prescribed in advance. The proof
goes by observing that the {\em positive series} (the series of
terms with positive sign) increases beyond any upper bound, and
that the {\em negative series} (the series of terms with negative
sign) decreases beyond any lower bound, when more and more terms
are added to the partial sums in the two series. Now, say that we
want to obtain the sum $S$ in a rearrangement of the alternating
harmonic series. Then we first take as many terms from the
positive series as are needed just to exceed $S$. Then take as
many of the terms from the negative series as needed just to come
below $S$.  Continue this way by taking as many terms as needed
from the positive series, from where we stopped earlier, just
again to exceed $S$. Then take as many terms as needed from the
negative series, from where we stopped earlier, just again to come
below $S$, and so on. Since $\frac 1n$ tends to $0$ as $n$
increases, the series just described is convergent with sum $S$.

The alternating harmonic series is the mathematical idea behind
the episode from {\em Alice's Adventures in Wonderland}, mentioned
at the beginning of this paper.

\section{Mathematics as a sixth sense} \label{section 9} \setzero

\vskip-5mm \hspace{5mm}

In 1982, {\em Dirac's string problem} was
presented on a shopping bag from a Danish supermarket chain
\cite{hans3}. This problem illustrates the property of half-spin
of certain elementary particles, mathematically predicted by the
physicist P.A.M. Dirac in the 1920s for elementary particles such
as the electron and the neutron. To convince colleagues sceptical
of his theory, Dirac conceived a model to illustrate a
corresponding phenomenon in the macroscopical world. This model
consists of a solid object (Dirac used a pair of scissors)
attached to two posts by loose (or elastic) strings, say with one
string from one end of the object and two strings from the other
end to the two posts. Dirac then demonstrated that a double
twisting of the strings could be removed by passing the strings
over and round the object, while he was not able to remove a
single twisting of the strings in this way. Rather deep
mathematics from topology is needed to explain the phenomenon, and
in this case mathematics comes in as a `sixth sense', by which we
`sense' (understand) counter-intuitive phenomena. An application
to the problem of transferring electrical current to a rotating
plate, without the wires getting tangled and breaking, was
patented in 1971.

\section{Making mathematics visible} \label{section 10} \setzero

\vskip-5mm \hspace{5mm}

In many countries there are increasing
demands to make the role and importance of the subjects taught in
the school system visible to the general public. If the curriculum
and the methods of teaching are not to stagnate, it is important
that an informed debate takes place in society about the
individual subjects in school.

For subjects strongly depending on mathematics, this is a great
challenge. Where mathematicians (and scientists in general) put
emphasis on explanations (proofs), the public (including
politicians) are mostly interested in results and consequences. It
is doubtless necessary to go part of the way to avoid technical
language in presenting mathematics for a wider audience.

On the initiative of the International Mathematical Union, the
year 2000 was declared {\em World Mathematical Year}. During the
year many efforts were made to reach the general public through
poster campaigns in metros, public lectures, articles in general
magazines, etc. and valuable experience was gained. The
mathematical year clearly demonstrated the value of an
international exchange of ideas in such matters, and the issue of
raising public awareness of mathematics is now on the agenda all
over the world from eight to infinity.

\label{lastpage}

\end{document}